\documentclass[preprint,review,12pt]{elsarticle}
\usepackage{amssymb}
\usepackage{amsfonts}
\usepackage{amsmath}
\usepackage{amsthm}
\usepackage{color}
\usepackage{graphicx}
\usepackage{epsfig,mathrsfs}
\usepackage[bf,SL,BF]{subfigure}
\usepackage{fancyhdr}
\usepackage{CJK}
\usepackage{caption}
\usepackage{wrapfig}
\usepackage{cases}

\usepackage{algorithm}
\usepackage{algorithmic}

\newtheorem*{lemma*}{Lemma A}

\numberwithin{equation}{section}

\renewcommand{\eqref}[1]{(\ref{#1})}

 \allowdisplaybreaks

\begin{document}
 \pagenumbering{arabic}

\begin{frontmatter}
\title{{\bf \noindent Second-order LOD multigrid method for multidimensional Riesz fractional diffusion equation}}
\author{Minghua Chen, Yantao Wang, Xiao Cheng, Weihua Deng$^{*}$}
\cortext[cor2]{Corresponding author. E-mail: dengwh@lzu.edu.cn.}
\address{School of Mathematics and Statistics,
Lanzhou University, Lanzhou 730000, P.R. China}

\date{}

\begin{abstract}
We propose a locally one dimensional (LOD) finite difference method for multidimensional Riesz fractional diffusion equation with variable coefficients on a finite domain. The numerical method is second-order convergent in both space and time directions, and its unconditional stability is strictly proved. Comparing with the popular first-order finite difference method for fractional operator, the form of obtained matrix algebraic equation is changed from $(I-A)u^{k+1}=u^k+b^{k+1}$ to $(I-{\widetilde A})u^{k+1}=(I+{\widetilde B})u^k+{\tilde b}^{k+1/2}$; the three matrices $A$, ${\widetilde A}$ and ${\widetilde B}$ are all Toeplitz-like, i.e.,  they have completely same structure and the computational count for matrix vector multiplication is $\mathcal{O}(N \mbox{log} N)$; and the computational costs for solving the two matrix algebraic equations are almost the same. The LOD-multigrid method is used to solve the resulting matrix algebraic equation, and the computational count is
$\mathcal{O}(N \mbox{log} N)$ and the required storage is $\mathcal{O}(N)$, where $N$ is the number of grid points. Finally, the extensive numerical experiments are performed to show the powerfulness of the second-order scheme and the LOD-multigrid method.





\end{abstract}
\end{frontmatter}

\section{Introduction}

In recent years considerable interests in fractional calculus have been
stimulated by the applications in physical, chemical, biological, and engineering, etc., areas \cite{Mainardi:10}. The definitions of fractional calculus are versatile, e.g., Riemann-Liouville derivative, Gr\"{u}nwald-Letnikov derivative, Caputo derivative, Weyl derivative and Riesz derivative et al \cite{Miller:93,Podlubny:99}, and they are not completely equivalent.  Depending on the particular applied field, sometimes one of its definitions is more popular than others. For example, the Riesz fractional derivative appears in the continuous limit of lattice models with long-range interactions \cite{Tarasov:10}. This paper focuses on the multidimensional Riesz fractional diffusion equations.

Nowadays, the finite difference discretization for space fractional derivatives is experiencing rapid development, including the Riesz fractional derivative; such as, Yang et al numerically study the Riesz space fractional PDEs with two different fractional orders $1<\alpha \le 2$ and $0 <\beta <1$ \cite{Yang:10}; Zhuang  et al consider a variable-order fractional advection-diffusion equation with a nonlinear source term on a finite domain \cite{Zhuang:09}. In the last two years, for the space fractional derivatives, we notice that two different second-order discretization schemes are developed \cite{Sousa:12, Tian:12}; even the third-order discrezation scheme is obtained \cite{Zhou:12} if a compact difference operator is performed on the discrezation scheme given in \cite{Tian:12}.

Another topic related to effectively solving the equations involving fractional operators is about how to efficiently solve the resulting matrix algebraic equations. The `unlucky' thing is that the matrix in the matrix algebraic equation is usually full because of the nonlocal properties of the fractional operators, and the `lucky' thing,  as pointed out in \cite{Pang:12, Wang:10, Wang:12}, is that the matrix has some special structure, i.e., the matrix is Toeplitz-like matrix, and the count of its matrix vector multiplication is $\mathcal{O}(N \mbox{log} N)$ by using the constructed circulant matrix and fast Fourier transform, and the required storage is $\mathcal{O}(N)$, where $N$ is the number of grid points. Pang and Sun \cite{Pang:12} successfully use the multigrid method (MGM) to efficiently solve the resulting matrix algebraic equation of the one dimensional fractional diffusion equation by using first-order discretization scheme \cite{Meerschaert:06}. Here we further extend the MGM to solve the matrix algebraic equation of the multidimensional Riesz fractional diffusion equation discretized by the second-order scheme.

We use the LOD strategy to solve the multidimensional Riesz fractional diffusion equation. The LOD methods include alternating direction (AD) methods and fractional step  procedures \cite{Douglas:01}. The AD methods were first introduced in three papers \cite{Douglas:55,DouPea:55,Peaceman:55}
 by Douglas, Peaceman, and Rachford.   The Peaceman and Rachford  (PR-AD) method  works well for two-dimensional problems.
 However, it can not be extended to higher dimensional problems. Douglas (D-AD) method  \cite{Douglas:55,Douglas:64,Douglas:01}
are valid for any dimensional equations. And PR-AD and D-AD are equivalent in two-dimensional problems. Both PR-AD and D-AD schemes are used in this paper to discretize the multidimensional Riesz fractional diffusion equation. In each dimension, the obtained matrix algebraic equation is solved by MGM. Although the spacial fractional derivative is discretized by second-order scheme, for any single dimension the form of the obtained matrix algebraic equation is $(I-{\widetilde A})u^{k+1}=(I+{\widetilde B})u^k+{\tilde b}^{k+1/2}$, in fact the corresponding form for the first-order discretization scheme is  $(I-A)u^{k+1}=u^k+b^{k+1}$; the three matrices $A$, ${\widetilde A}$ and ${\widetilde B}$ are all Toeplitz-like, and they have completely same structure and the computational count for matrix vector multiplication is $\mathcal{O}(N \mbox{log} N)$;  and the computational costs for solving the two matrix algebraic equations are almost the same. In other words, the second-order scheme improves the accuracy but almost without increasing the computational cost.

%

More concretely, in this paper using the second-order accurate and unconditionally stable computational scheme and LOD-MGM, we solve the following multidimensional variable coefficients Riesz fractional diffusion equation
with the computational count $\mathcal{O}(N \mbox{log} N)$ and storage $\mathcal{O}(N)$,
\begin{equation} \label{1.1}
\left\{ \begin{split}
\frac{\partial u(x,y,z,t) }{\partial t}&=c(x,y,z,t)   \frac{\partial ^{\alpha}u(x,y,z,t)}{\partial |x|^{\alpha}}
                                                           +d(x,y,z,t)\frac{\partial ^{\beta}u(x,y,z,t)}{\partial |y|^{\beta}}\\
                                         &\quad +e(x,y,z,t)\frac{\partial ^{\gamma}u(x,y,z,t)}{\partial |z|^{\gamma}}+f(x,y,z,t),\\
                            u(x,y,z,0) &=u_0(x,y,z) ~~~~\, {\rm for}~~~ (x,y,z) \in \Omega,\\
                             u(x,y,z,t)&=0 ~~~ {\rm for}~~~ (x,y,z,t) \in \partial \Omega \times (0,T],
 \end{split}
 \right.
\end{equation}
in the domain $\Omega=(x_L,x_R) \times (y_L,y_R)\times (z_L,z_R),\, 0< t \leq T$,
where the orders of the Riesz fractional derivatives are  $1<\alpha,\beta,\gamma< 2$;
$f(x,y,z,t)$ is a source term and the variable coefficients
$c(x,y,z,t)\geq 0$, $d(x,y,z,t)\geq0$, $e(x,y,z,t)\geq0$; the Riesz fractional derivative  for $n \in \mathbb{N}$, $n-1 < \nu \leq n$,
is defined as  \cite{Chen:12,Tarasov:10}
\begin{equation}\label{1.2}
\frac{\partial ^{\nu}u(x,y,z,t)}{\partial |x|^{\nu}} =-\kappa_{\nu}\left( _{x_L}\!D_x^{\nu}+ _{x}\!D_{x_R}^{\nu} \right)u(x,y,z,t),
\end{equation}
where the coefficient $\kappa_{\nu}=\frac{1}{2cos(\nu \pi/2)}$, and
\begin{equation}\label{1.3}
 _{x_L}D_x^{\nu}u(x,y,z,t)=
\frac{1}{\Gamma(n-\nu)} \displaystyle \frac{\partial^n}{\partial x^n}
 \int_{x_L}\nolimits^x{\left(x-\xi\right)^{n-\nu-1}}{u(\xi,y,z,t)}d\xi,
\end{equation}

\begin{equation}\label{1.4}
 _{x}D_{x_R}^{\nu}u(x,y,z,t)=
 \frac{(-1)^n}{\Gamma(n-\nu)}\frac{\partial^n}{\partial x^n}
\int_{x}\nolimits^{x_R}{\left(\xi-x\right)^{n-\nu-1}}{u(\xi,y,z,t)}d\xi,
\end{equation}
are the left and right Riemann-Liouville space fractional
derivatives, respectively.
%

The  outline of this paper is as follows. In the next section, we introduce the
 second-order finite difference discretizations for the
Riesz fractional derivatives; and the full discretization of (\ref{1.1}) is
derived, where the Crank-Nicolson scheme and LOD method are combined together. We theoretically prove the presented finite
difference scheme is unconditionally stable in Section 3.
In Section 4 we propose a V-cycle LOD-MGM for the resulting system of (\ref{1.1}).
 To show the powerfulness of the second-order scheme and LOD-MGM, the extensive numerical experiments are performed in Section 5.
Finally, we conclude the paper with some remarks in the last section.

\section{Derivation of the finite difference scheme}\label{sec:1}
In this section, we derive the full discretization schemes of (\ref{1.1}).
The first subsection introduces   the
 second-order finite difference discretizations for the
Riesz fractional    derivatives  in a finite domain.
Then in the second subsection, we present the scheme for the one dimensional case of
(\ref{1.1}). The third and  fourth  subsections detailedly provide the
 two dimensional case of (\ref{1.1}) and (\ref{1.1}) itself, respectively.

\subsection{Discretizations for the  Riesz fractional derivatives}
Take the mesh points $x_i=x_L+i\Delta x,i=0,1,\ldots ,{N_x}$, $y_j=y_L+j\Delta y,j=0,1,\ldots ,{N_y}$, $z_l=z_L+l\Delta z,l=0,1,\ldots ,{N_z}$
and $t_k=k\Delta t,k=0,1,\ldots ,{N_t}$, where
 $\Delta x=(x_R-x_L)/{N_x}$, $\Delta y=(y_R-y_L)/{N_y}$, $\Delta z=(z_R-z_L)/{N_z}$, $\Delta t=T/{N_t}$,
 i.e., $\Delta x$, $\Delta y$ and $\Delta z$ are the uniform space stepsizes in the corresponding directions, $\Delta t$ the time stepsize.
 For $\nu\in (1,2)$,  the left and right  Riemann-Liouville space fractional
derivatives (\ref{1.3}) and (\ref{1.4}) have the second-order approximation operators
$\delta_{\nu,_+x}{u_{i,j,l}^k}$ and $\delta_{\nu,_-x}{u_{i,j,l}^k}$, respectively, given in a finite
domain \cite{Chen:12,Sousa:12}, where ${u_{i,j,l}^k}$ denotes the
approximated value of $u(x_i,y_j,z_l,t_k)$.

The approximation operator of (\ref{1.3}) is defined by \cite{Chen:12,Sousa:12}
\begin{equation}\label{2.1}
\begin{split}
  \delta_{\nu,_+x}{u_{i,j,l}^k}:=\frac{1}{\Gamma(4-\nu)(\Delta x)^{\nu}}\sum_{m=0}^{i+1} g_m^{\nu}u_{i-m+1,j,l}^k,
  \end{split}
\end{equation}
and there exists
\begin{equation}\label{2.2}
  _{x_L}\!D_x^{\nu}u(x,y,z,t)= \delta_{\nu,_+x}{u_{i,j,l}^k}+\mathcal{O}(\Delta x)^2,
\end{equation}
where
\begin{equation}\label{2.3}
g_m^{\nu}=\left\{ \begin{array}
 {l@{\quad } l}

1,&m =0,\\

-4+2^{3-\nu},&m =1,\\

6-2^{5-\nu}+3^{3-\nu},&m =2,\\

(m+1)^{3-\nu}-4m^{3-\nu}+6(m-1)^{3-\nu}\\
~~~~~~~~~~~-4(m-2)^{3-\nu}+(m-3)^{3-\nu},&m \geq 3.\\
 \end{array}
 \right.
\end{equation}

Analogously,   the   approximation operator of (\ref{1.4}) is  described as \cite{Chen:12}
 \begin{equation}\label{2.4}
  \delta_{\nu,_-x}{u_{i,j,l}^k}:
                              =\frac{1}{\Gamma(4-\nu)(\Delta x)^{\nu}}\sum_{m=0}^{N_x-i+1}g_m^{\nu}u_{i+m-1,j,l}^k,
\end{equation}
and it holds that
\begin{equation}\label{2.5}
  _{x}D_{x_R}^{\nu}u(x,y,z,t)= \delta_{\nu,_-x}{u_{i,j,l}^k}+\mathcal{O}(\Delta x)^2,
\end{equation}
where $g_m^{\nu}$ is defined by (\ref{2.3}).

Combining (\ref{2.2}) and (\ref{2.5}), we obtain the approximation operator of the (Riemann-Liouville) Riesz fractional derivative

\begin{equation}\label{2.6}
\begin{split}
\frac{\partial ^{\nu}u(x,y_j,z_l,t_k)}{\partial |x|^{\nu}}\Big|_{x=x_i}
 &=-\kappa_{\nu}\left( _{x_L}\!D_x^{\nu}+ _{x}\!D_{x_R}^{\nu} \right)u(x,y_j,z_l,t_k)\big|_{x=x_i}\\
& =
 -\kappa_{\nu} \left( \delta_{\nu,_+x}+\delta_{\nu,_-x} \right) {u_{i,j,l}^k} +\mathcal{O}(\Delta x)^2\\
&=\frac{-\kappa_{\nu}}{\Gamma(4-\nu)\Delta x^{\nu}}
  \left(\sum_{m=0}^{i+1}g_m^{\nu}u_{i-m+1,j,l}^k+\sum_{m=0}^{N_x-i+1}g_m^{\nu}u_{i+m-1,j,l}^k \right)+\mathcal{O}(\Delta x)^2\\
&=\frac{-\kappa_{\nu}}{\Gamma(4-\nu)\Delta x^{\nu}}
  \left(\sum_{m=0}^{i+1}g_{i-m+1}^{\nu}u_{m,j,l}^k+\sum_{m=i-1}^{N_x}g_{m-i+1}^{\nu}u_{m,j,l}^k \right)+\mathcal{O}(\Delta x)^2\\
  &:=\frac{-\kappa_{\nu}}{\Gamma(4-\nu)\Delta x^{\nu}}
  \sum_{m=0}^{N_x}\widetilde{g}_{i,m}^{\nu}u_{m,j,l}^k+\mathcal{O}(\Delta x)^2,\\
\end{split}
\end{equation}
where
\begin{equation}\label{2.7}
\widetilde{g}_{i,m}^{\nu}=\left\{ \begin{array}
 {l@{\quad } l}
  g_{i-m+1}^{\nu},&m < i-1,\\
  g_{0}^{\nu}+g_{2}^{\nu} ,&m=i-1,\\
 2g_{1}^{\nu},&m=i,\\
 g_{0}^{\nu}+g_{2}^{\nu} ,&m=i+1,\\
 g_{m-i+1}^{\nu} ,&m>i+1,\\
 \end{array}
 \right.
\end{equation}
with $i=1,\ldots,N_x-1$, together with the Dirichlet boundary conditions that define $u_{0,j,l}^k$ and $u_{N_x,j,l}^k$  as appropriate.

Taking  $\nu=2$, both  Eq. (\ref{2.2}) and (\ref{2.5}) reduce to the following form
\begin{equation*}
 \frac{\partial^2 u(x_i,y,z,t)}{\partial x^2}=\frac{u(x_{i+1},y,z,t)-2u(x_i,y,z,t)+u(x_{i-1},y,z,t)}{(\Delta x)^2}+\mathcal{O}(\Delta x)^2.
\end{equation*}

Similarly, it is easy to get the one-dimensional and two-dimensioanl case of (\ref{2.1})-(\ref{2.7}).

\subsection{Numerical scheme for 1D}

Consider  the one-dimensional Riesz fractional   diffusion equation
\begin{equation}\label{2.8}
\frac{\partial u(x,t) }{\partial t}=c(x,t)\frac{\partial ^{\alpha}u(x,t)}{\partial |x|^{\alpha}}+f(x,t).
\end{equation}

In the time direction, we use the Crank-Nicolson scheme.
Taking the uniform time step $\Delta t$ and space step $\Delta x$,
and setting $c_i^k=c(x_i,t_k)$ and
$f_i^{k+1/2}=f(x_i,t_{k+1/2})$, where $t_{k+1/2}=(t_k+t_{k+1})/2$,
the full discretization of (\ref{2.8}) has the following form
\begin{equation}\label{2.9}
\begin{split}
\frac{u_i^{k+1}-u_i^k}{\Delta t}=\frac{-\kappa_{\alpha}c_i^{k+1/2}}{\Gamma(4-\alpha)\Delta x^{\alpha}}
  \sum_{m=0}^{N_x}\widetilde{g}_{i,m}^{\alpha}\frac{u_m^k +u_m^{k+1}}{2} +f_i^{k+1/2}.
\end{split}
\end{equation}
 Then  (\ref{2.9}) can be expressed as
\begin{equation}\label{2.10}
 \left (1-\frac{\Delta t}{2} \delta''_{\alpha,_x}     \right )u_i^{k+1}=
 \left (1+\frac{\Delta t}{2}  \delta''_{\alpha,_x}  \right )u_i^{k}+ \Delta t f_i^{k+1/2},
\end{equation}
where
\begin{equation*}
 \delta''_{\alpha,_x}{u_i^k}:=\frac{-\kappa_{\alpha}c_i^{k+1/2}}{\Gamma(4-\alpha)\Delta x^{\alpha}}\sum_{m=0}^{N_x}\widetilde{g}_{i,m}^{\alpha}u_m^k;
\quad      \delta''_{\alpha,_x}{u_i^{k+1}}:=\frac{-\kappa_{\alpha}c_i^{k+1/2}}{\Gamma(4-\alpha)\Delta x^{\alpha}}\sum_{m=0}^{N_x}\widetilde{g}_{i,m}^{\alpha}u_m^{k+1}.
\end{equation*}

  Putting $\xi_i^{k+1/2}=\frac{-\Delta t\kappa_{\alpha}c_i^{k+1/2}}{2\Gamma(4-\alpha)\Delta x^{\alpha}}$,
 the system of equations given by (\ref{2.10}) takes the form
\begin{equation}\label{2.11}
   (I-A^{k+1/2})U^{k+1}=(I+A^{k+1/2})U^{k}+\Delta t F^{k+1/2},
\end{equation}
 where $I$ is the $(N_x-1) \times (N_x-1) $ identity  matrix,
 \begin{equation*}
 U^{k}=[u_1^k,u_2^k,\ldots,u_{N_x-1}^k]^{\rm T}, ~~F^{k+1/2}=[f_1^{k+1/2},f_2^{k+1/2},\ldots,f_{N_x-1}^{k+1/2}]^{\rm T},
  \end{equation*}
 and the discretizations at the interior $x$-gridpoints define the entries of the matrix $A^{k+1/2}$,
 $A_{i,m}^{k+1/2}$ for $i=1,\ldots,N_x-1$ and $m=1,\ldots,N_x-1 $
  are defined by
  \begin{equation}\label{2.12}
A_{i,m}^{k+1/2}=\left\{ \begin{array}
 {l@{\quad} l}
g_{i-m+1}^{\alpha}  \xi_i^{k+1/2},&m<i-1,\\
(g_0^{\alpha}+g_2^{\alpha})  \xi_i^{k+1/2},&m=i-1,\\
 2g_1^{\alpha} \xi_i^{k+1/2},& m=i,\\
(g_0^{\alpha}+g_2^{\alpha}) \xi_i^{k+1/2},&m=i+1,\\
g_{m-i+1}^{\alpha} \xi_i^{k+1/2},&m>i+1.\\
 \end{array}
 \right.
\end{equation}

\subsection{LOD scheme for 2D}

Consider the following two-dimensional Riesz fractional diffusion equation
\begin{equation}\label{2.13}
\frac{\partial u(x,y,t) }{\partial t}=c(x,y,t)   \frac{\partial ^{\alpha}u(x,y,t)}{\partial |x|^{\alpha}}
+d(x,y,t)   \frac{\partial ^{\beta}u(x,y,t)}{\partial |y|^{\beta}} +f(x,y,t).
\end{equation}

 Analogously we still use the Crank-Nicolson scheme to do the discretization in time direction.
 Taking $u_{i,j}^k$ as the approximated value of $u(x_i,y_j,t_k)$,
 $c_{i,j}^k=c(x_i,y_j,t_k)$,  $d_{i,j}^k=d(x_i,y_j,t_k)$, $t_{n+1/2}=(t_n+t_{n+1})/2$, $f_{i,j}^{k+1/2}=f(x_i,y_j,t_{k+1/2})$,
 $\Delta x=(x_R-x_L)/{N_x}$, and $\Delta y=(y_R-y_L)/{N_y}$, for the uniform space steps $\Delta x,\Delta y$
 and the time stepsize $\Delta t$, the resulting discretization of (\ref{2.13}) can be written as
\begin{equation}\label{2.14}
\begin{split}
  \frac{u_{i,j}^{k+1}-u_{i,j}^k}{\Delta t}
  =&\frac{-\kappa_{\alpha}c_{i,j}^{k+1/2}}{\Gamma(4-\alpha)\Delta x^{\alpha}}
  \sum_{m=0}^{N_x}\widetilde{g}_{i,m}^{\alpha}\frac{u_{m,j}^k +u_{m,j}^{k+1}}{2}\\
  &+\frac{-\kappa_{\beta}d_{i,j}^{k+1/2}}{\Gamma(4-\beta)\Delta x^{\beta}}
  \sum_{m=0}^{N_y} \widetilde{g}_{j,m}^{\beta}\frac{u_{i,m}^k +u_{i,m}^{k+1}}{2}+f_{i,j}^{k+1/2}.
\end{split}
\end{equation}

Similarly,  we define
\begin{equation*}
\begin{split}
&\delta''_{\alpha,_x}{u_{i,j}^k}=\frac{-\kappa_{\alpha}c_{i,j}^{k+1/2}}{\Gamma(4-\alpha)\Delta x^{\alpha}}
                                         \sum_{m=0}^{N_x}\widetilde{g}_{i,m}^{\alpha}u_{m,j}^k;
\quad\delta''_{\alpha,_x}{u_{i,j}^{k+1}}=\frac{-\kappa_{\alpha}c_{i,j}^{k+1/2}}{\Gamma(4-\alpha)\Delta x^{\alpha}}
                                         \sum_{m=0}^{N_x}\widetilde{g}_{i,m}^{\alpha}u_{m,j}^{k+1};\\
&\delta''_{\beta,_y}{u_{i,j}^k}=\frac{-\kappa_{\beta}d_{i,j}^{k+1/2}}{\Gamma(4-\beta)\Delta y^{\beta}}
                                         \sum_{m=0}^{N_y}\widetilde{g}_{j,m}^{\beta}u_{i,m}^k;
\quad\delta''_{\beta,_y}{u_{i,j}^{k+1}}=\frac{-\kappa_{\beta}d_{i,j}^{k+1/2}}{\Gamma(4-\beta)\Delta y^{\beta}}
                                         \sum_{m=0}^{N_y}\widetilde{g}_{j,m}^{\beta}u_{i,m}^{k+1};\\
\end{split}
\end{equation*}
then Eq. (\ref{2.14}) can be rewritten  as
\begin{equation}\label{2.15}
\begin{split}
 &\left (1-\frac{\Delta t}{2}\delta''_{\alpha,_x} -\frac{\Delta t}{2}\delta''_{\beta,_y} \right )u_{i,j}^{k+1}
 =\left (1+\frac{\Delta t}{2}\delta''_{\alpha,_x} +\frac{\Delta t}{2}\delta''_{\beta,_y} \right )u_{i,j}^{k}+ \Delta t f_{i,j}^{k+1/2}.
\end{split}
\end{equation}

For the two-dimensional Riesz fractional  diffusion
equation (\ref{2.13}), the relevant perturbation of (\ref{2.15}) is
of the form
\begin{equation}\label{2.16}
\begin{split}
 &\left (1-\frac{\Delta t}{2}\delta''_{\alpha,x} \right )\left(1-\frac{\Delta t}{2}\delta''_{\beta,y}\right )u_{i,j}^{k+1}
=\left (1+\frac{\Delta t}{2}\delta''_{\alpha,x}\right)\left(1+\frac{\Delta t}{2}\delta''_{\beta,y}\right )u_{i,j}^{k}+ \Delta t f_{i,j}^{k+1/2}. \\
\end{split}
\end{equation}

The scheme (\ref{2.16}) differs from (\ref{2.15}) by a perturbation \cite{Tadjeran:07}
\begin{equation*}
  \frac{(\Delta
  t)^2}{4}\delta''_{\alpha,x}\delta''_{\beta,y}(u_{i,j}^{k+1}-u_{i,j}^k),
\end{equation*}
which may be deduced by distributing the operator products in
(\ref{2.16}). Since $(u_{i,j}^{k+1}-u_{i,j}^k)$ is an
$\mathcal{O}(\Delta t)$ term, it follows that the perturbation
contributes an $\mathcal{O}((\Delta t)^2)$ error component to the
truncation error of (\ref{2.15}). Thus, the scheme (\ref{2.16}) has
a truncation error also $\mathcal{O}((\Delta x)^2)+\mathcal{O}((\Delta y)^2)+\mathcal{O}((\Delta t)^2)$.

For efficiently solving system (\ref{2.16}), the following techniques can be used:

 D-AD    scheme \cite{Douglas:55,Douglas:01}:
\begin{align}
\left (1-\frac{\Delta t}{2}\delta''_{\alpha,x} \right )u_{i,j}^{*}
                &=\left(1+\frac{\Delta t}{2}\delta''_{\alpha,x}+\Delta t\delta''_{\beta,y}\right )u_{i,j}^{k}+\Delta t f_{i,j}^{k+1/2}; \label{2.17}\\
\left(1-\frac{\Delta t}{2}\delta''_{\beta,y}\right )u_{i,j}^{k+1}
                &=u_{i,j}^{*}- \frac{\Delta t}{2}\delta''_{\beta,y}u_{i,j}^{k} \label{2.18};
\end{align}

where $u_{i,j}^{*}$ is an intermediate solution. Subtracting (\ref{2.18}) from  ({\ref{2.17}), we obtain
\begin{equation*}
u_{i,j}^{*}= u_{i,j}^{k+1}+\frac{\Delta t}{2}\delta''_{\beta,y}\left(u_{i,j}^{k}-u_{i,j}^{k+1}\right).
\end{equation*}

PR-AD scheme \cite{Peaceman:55}:

\begin{align}
\left (1-\frac{\Delta t}{2}\delta''_{\alpha,x} \right )u_{i,j}^{*}
 &=\left(1+\frac{\Delta t}{2}\delta''_{\beta,y}\right )u_{i,j}^{k}+ \frac{\Delta t}{2}f_{i,j}^{k+1/2}; \label{2.19}\\
\left(1-\frac{\Delta t}{2}\delta''_{\beta,y}\right )u_{i,j}^{k+1}
&=\left (1+\frac{\Delta t}{2}\delta''_{\alpha,x}\right)u_{i,j}^{*}+\frac{\Delta t}{2} f_{i,j}^{k+1/2}; \label{2.20}
\end{align}
with  intermediate solution $u_{i,j}^{*}$.
Subtracting (\ref{2.20}) from  ({\ref{2.19}), we have
\begin{equation*}
2u_{i,j}^{*}=\left (1-\frac{\Delta t}{2}\delta''_{\beta,y} \right)u_{i,j}^{k+1}
 +\left(1+\frac{\Delta t}{2}\delta''_{\beta,y}\right )u_{i,j}^{k}.
\end{equation*}

\subsection{D-AD scheme for 3D}

 Similarly, the resulting discretization of (\ref{1.1}) can be written as,
\begin{equation}\label{2.21}
\begin{split}
 &\left (1-\frac{\Delta t}{2}\delta''_{\alpha,x}-\frac{\Delta t}{2}\delta''_{\beta,y}-\frac{\Delta t}{2}\delta''_{\gamma,z}\right )u_{i,j,l}^{k+1}\\
 &\quad=\left (1+\frac{\Delta t}{2}\delta''_{\alpha,x}+\frac{\Delta t}{2}\delta''_{\beta,y}+\frac{\Delta t}{2}\delta''_{\gamma,z}\right )u_{i,j,l}^{k}
     +f_{i,j,l}^{k+1/2}\Delta t.
\end{split}
\end{equation}
The perturbation equation of (\ref{2.21}) is of the form
\begin{equation}\label{2.22}
\begin{split}
 &  \left(1-\frac{\Delta t}{2}\delta''_{\alpha,x}\right)      \left(1-\frac{\Delta t}{2}\delta''_{\beta,y}\right)
    \left(1-\frac{\Delta t}{2} \delta''_{\gamma,z}\right)  u_{i,j,l}^{k+1}\\
 &\quad =  \left (1+\frac{\Delta t}{2}\delta''_{\alpha,x}\right) \left(1+\frac{\Delta t}{2}\delta''_{\beta,y}\right )
          \left(1+\frac{\Delta t}{2}\delta''_{\gamma,z}\right) u_{i,j,l}^{k}
 +f_{i,j,l}^{k+1/2}\Delta t.
\end{split}
\end{equation}
The scheme (\ref{2.22}) differs from (\ref{2.21}) by the perturbation term
$$ \frac{(\Delta
  t)^2}{4}(\delta''_{\alpha,x}\delta''_{\beta,y}+\delta''_{\alpha,x}\delta''_{\gamma,z}+\delta''_{\beta,y}\delta''_{\gamma,z})(u_{i,j,l}^{k+1}-u_{i,j,l}^k) - \frac{(\Delta
  t)^3}{8}\delta''_{\alpha,x}\delta''_{\beta,y}\delta''_{\gamma,z}(u_{i,j,l}^{k+1}-u_{i,j,l}^k).$$

The system of the equations defined by (\ref{2.22}) can be solved by the
 D-AD  scheme \cite{Douglas:64,Douglas:01}
\begin{equation}\label{2.23}
\left (1-\frac{\Delta t}{2}\delta''_{\alpha,x} \right )u_{i,j,l}^{k,1}
=\left(1+\frac{\Delta t}{2}\delta''_{\alpha,x}+\Delta t\delta''_{\beta,y}+\Delta t\delta''_{\gamma,z}\right )u_{i,j,l}^{k}+\Delta t f_{i,j,l}^{k+1/2};
\end{equation}
\begin{equation}\label{2.24}
\left(1-\frac{\Delta t}{2}\delta''_{\beta,y}\right )u_{i,j,l}^{k,2}=u_{i,j,l}^{k,1}- \frac{\Delta t}{2}\delta''_{\beta,y}u_{i,j,l}^{k};
\end{equation}
\begin{equation}\label{2.25}
\left(1-\frac{\Delta t}{2}\delta''_{\gamma,z}\right )u_{i,j,l}^{k+1}=u_{i,j,l}^{k,2}- \frac{\Delta t}{2}\delta''_{\gamma,z}u_{i,j,l}^{k}.
\end{equation}

For maintaining the consistency, we need  to carefully specify the boundary conditions of $u_{i,j,l}^{n,1}$ and $u_{i,j,l}^{k,2}$.
According to  ({\ref{2.23})-({\ref{2.25}), we obtain
 \begin{equation*}
 \begin{split}
&u_{i,j,l}^{k,1}= u_{i,j,l}^{k+1}+\frac{\Delta t}{2}\left(\delta''_{\beta,y}+\delta''_{\gamma,z}\right)\left(u_{i,j,l}^{k}-u_{i,j,l}^{k+1}\right)
+\frac{(\Delta t)^2}{4}\delta''_{\beta,y}\delta''_{\gamma,z}u_{i,j,l}^{k+1};\\
&u_{i,j,l}^{k,2}= u_{i,j,l}^{k+1}+\frac{\Delta t}{2}\delta''_{\gamma,z}\left(u_{i,j,l}^{k}-u_{i,j,l}^{k+1}\right).\\
\end{split}
\end{equation*}

\section{Convergence and Stability Analysis}
We show the convergence for one-dimensional and multidimensional
Riesz fractional  diffusion equation by proving the
consistency and stability (according to Lax's equivalence theorem).

\noindent{\bf Lemma 3.1}\label{Lemma:1} (\cite{Chen:12}).
\emph{The coefficients $\widetilde{g}_{i,m}^{\nu}$, $\nu \in(1,2)$ defined in (\ref{2.7}) satisfy}

\begin{equation*}\label{3.27}
  \begin{split}
&(1) ~~ \widetilde{g}_{i,i}^{\nu}<0,\,\,\,\,\,\, \widetilde{g}_{i,m}^{\nu} > 0\,\,(m\neq i);\\
&(2) ~~ \sum\limits_{m=0}^{N_x}\widetilde{g}_{i,m}^{\nu}< 0 ~~ \mbox {and} ~~  -\widetilde{g}_{i,i}^{\nu}>\!\!\!\!\!
         \sum\limits_{m=0,m\neq i}^{N_x}\!\!\!\!\widetilde{g}_{i,m}^{\nu}.
  \end{split}
\end{equation*}

\subsection{The stability of the numerical methods in 1D}

\noindent{\bf Theorem 3.2.}\label{Theorem:2}
The Crank-Nicholson  scheme (\ref{2.11}) of the Riesz fractional  diffusion
equation (\ref{2.9}) with $1<\alpha<2$ is unconditionally stable.

\begin{proof}
First, we prove that the eigenvalues of the matrix $A^{k+1/2}$ have negative real parts.
Note that $A_{i,i}^{k+1/2}= \widetilde{g}_{i,i}^{\alpha} \xi_i^{k+1/2}$, and  from Lemma 3.1 we obtain
\begin{equation}
  r_i= \sum\limits_{m=0,m\neq i}^{N_x}\!\!\!\! |A_{i,m}^{k+1/2}|=\xi_i^{k+1/2}\!\!\!\! \sum\limits_{m=0,m\neq i}^{N_x}\!\!\!\!\widetilde{g}_{i,m}^{\alpha}
      < -A_{i,i}^{k+1/2}.
\end{equation}
According to the Greschgorin theorem \cite{Isaacson:66}, the eigenvalues of the matrix $A^{k+1/2}$ are in the disks centered at $A_{i,i}^{k+1/2}$, with radius
$r_i$, i.e.,  the eigenvalues $\lambda$ of the matrix  $A^{k+1/2}$ satisfy
\begin{equation}
  |\lambda -A_{i,i}^{k+1/2} | \leq r_i,
\end{equation}
thus,  the eigenvalues of the matrix $A^{k+1/2}$ have negative real parts.
Similarly, we can prove that the eigenvalues of the matrix $I-A^{k+1/2}$ have a magnitude greater than $1$ and invertible.

Note that
 $\lambda$ is an eigenvalue of the matrix $A^{k+1/2}$ if and only if
$1-\lambda$ is an eigenvalue of the matrix $I-A^{k+1/2}$, if and only if $(1-\lambda)^{-1}(1+\lambda)$ is an
eigenvalue of the matrix $(I-A^{k+1/2})^{-1}(I+A^{k+1/2})$. Since  $\Re(\lambda) <0 $, it implies that $|(1-\lambda)^{-1}(1+\lambda)|<1$.
Hence, the spectral radius of the matrix $(I-A^{k+1/2})^{-1}(I+A^{k+1/2})$ is less than $1$.
\end{proof}

\subsection{The stability of the numerical methods in 2D}
Under a commutativity assumption for the operators $\left (1-\frac{\Delta t}{2}\delta''_{\alpha,x} \right )$
and $\left(1-\frac{\Delta t}{2}\delta''_{\beta,y}\right )$ in (\ref{2.15}),
the PR-AD scheme and  D-AD scheme  will
be shown to be unconditionally stable. The commutativity assumption for these two operators is a common practice in establishing
stability of the classical AD methods for the diffusion \cite{Douglas:01,Tadjeran:07}.
 The commutativity of these operators implies that the matrices $A_{2D,{\Delta x}}^{k+1/2}$
and $A_{2D,{\Delta y}}^{k+1/2}$ given in (\ref{3.7}) commute.

\noindent{\bf Theorem 3.3.}\label{Theorem:3}
Both the D-AD scheme (\ref{2.17})-(\ref{2.18}) and PR-AD scheme (\ref{2.19})-(\ref{2.20}),
defined by (\ref{2.16}),  are unconditionally stable for   $\alpha,\beta \in (1,2)$,
 if the matrices $A_{2D,{\Delta x}}^{k+1/2}$ and $A_{2D,{\Delta y}}^{k+1/2}$  commute.

\begin{proof}
D-AD scheme (\ref{2.17})-(\ref{2.18}) can be expressed in the form
\begin{align}
&(I-A_{2D,{\Delta x}}^{k+1/2})U^{*}=(I+A_{2D,{\Delta x}}^{k+1/2}+2A_{2D,{\Delta y}}^{k+1/2})U^{k}+\Delta t F^{k+1/2}; \label{3.3}\\
&(I-A_{2D,{\Delta y}}^{k+1/2})U^{k+1}=U^{*}  - A_{2D,{\Delta y}}^{k+1/2}U^{k}; \label{3.4}
\end{align}
and PR-AD scheme (\ref{2.19})-(\ref{2.20}) is of the form
\begin{align}
&(I-A_{2D,{\Delta x}}^{k+1/2})U^{*}=(I+A_{2D,{\Delta y}}^{k+1/2})U^{k}+\frac{\Delta t}{2} F^{k+1/2}; \label{3.5}\\
&(I-A_{2D,{\Delta y}}^{k+1/2})U^{k+1}=(I+A_{2D,{\Delta x}}^{k+1/2})U^{*}+\frac{\Delta t}{2} F^{k+1/2}; \label{3.6}
\end{align}
where the matrices $A_{2D,{\Delta x}}^{k+1/2}$ and $A_{2D,{\Delta y}}^{k+1/2}$ denote the
operators $\frac{\Delta t}{2}\delta''_{\alpha,x}$ and $\frac{\Delta t}{2}\delta''_{\beta,y} $, and
\begin{equation*}
\begin{split}
U^{k}=[u_{1,1}^k,u_{2,1}^k,\ldots,u_{N_x-1,1}^k,u_{1,2}^k,u_{2,2}^k,\dots,u_{N_x-1,2}^k,\ldots,u_{1,N_y-1}^k,u_{2,N_y-1}^k,\ldots,u_{N_x-1,N_y-1}^k]^T,\\
U^{*}=[u_{1,1}^*,u_{2,1}^*,\ldots,u_{N_x-1,1}^*,u_{1,2}^*,u_{2,2}^*,\dots,u_{N_x-1,2}^*,\ldots,u_{1,N_y-1}^*,u_{2,N_y-1}^*,\ldots,u_{N_x-1,N_y-1}^*]^T,
\end{split}
\end{equation*}
and the vector $F^{k+1/2}$ absorbs the source terms $f_{i,j}^{k+1/2}$ and the Dirichlet boundary conditions at time $t=t_{k+1}$ in the discretized equation. The matrices $A_{2D,{\Delta x}}^{k+1/2}$ and $A_{2D,{\Delta y}}^{k+1/2}$ are matrices of size $(N_x-1)(N_y-1)\times(N_x-1)(N_y-1)$.

Let us cancel the intermediate solution $U^*$, then D-AD scheme and PR-AD scheme  have the same form
  \begin{equation}\label{3.7}
   (I-A_{2D,{\Delta x}}^{k+1/2})(I-A_{2D,{\Delta y}}^{k+1/2})U^{k+1}=(I+A_{2D,{\Delta x}}^{k+1/2})(I+A_{2D,{\Delta y}}^{k+1/2})U^{k}+\Delta tF^{k+1/2},
  \end{equation}
from  (\ref{3.7}) we have the perturbation equation
\begin{equation*}
   (I-A_{2D,{\Delta x}}^{k+1/2})(I-A_{2D,{\Delta y}}^{k+1/2})E^{k+1}=(I+A_{2D,{\Delta x}}^{k+1/2})(I+A_{2D,{\Delta y}}^{k+1/2})E^{k},
  \end{equation*}
where
\begin{equation*}
E^{k}=[e_{1,1}^k,e_{2,1}^k,\ldots,e_{N_x-1,1}^k,e_{1,2}^k,e_{2,2}^k,\dots,
e_{N_x-1,2}^k,\ldots,e_{1,N_y-1}^k,e_{2,N_y-1}^k,\ldots,e_{N_x-1,N_y-1}^k]^T,
\end{equation*}
and  $e_{i,j}^k=u(x_i,y_j,t_k)-u_{i,j}^k$,  consequently
\begin{equation*}
    E^k=[(I-A_{2D,{\Delta y}}^{k+1/2})^{-1}(I-A_{2D,{\Delta x}}^{k+1/2})^{-1}(I+A_{2D,{\Delta x}}^{k+1/2})(I+A_{2D,{\Delta y}}^{k+1/2})]^kE^0.
  \end{equation*}
 Since the matrices $A_{2D,{\Delta x}}^{k+1/2}$ and $A_{2D,{\Delta y}}^{k+1/2}$ commute, it can be written as
    \begin{equation*}
    E^k=[(I-A_{2D,{\Delta x}}^{k+1/2})^{-1}(I+A_{2D,{\Delta x}}^{k+1/2})]^k[(I-A_{2D,{\Delta y}}^{k+1/2})^{-1}(I+A_{2D,{\Delta y}}^{k+1/2})]^kE^0.
  \end{equation*}
 According to  Theorem 3.2,  similarly, it is easy to check that
 the eigenvalues of the matrix $A_{2D,{\Delta x}}^{k+1/2}$ and $A_{2D,{\Delta y}}^{k+1/2}$ have negative real parts.
 Then,  both the spectral radius of the matrixes $(I-A_{2D,{\Delta x}}^{k+1/2})^{-1}(I+A_{2D,{\Delta x}}^{k+1/2})$
  and $(I-A_{2D,{\Delta y}}^{k+1/2})^{-1}(I+A_{2D,{\Delta y}}^{k+1/2})$ are less than 1,
therefore the sequence
  $[(I-A_{2D,{\Delta x}}^{k+1/2})^{-1}(I+A_{2D,{\Delta x}}^{k+1/2})]^k$
  and $[(I-A_{2D,{\Delta y}}^{k+1/2})^{-1}(I+A_{2D,{\Delta y}}^{k+1/2})]^k$  converge to zero matrix \cite{Saad:03}.
   Hence, the difference  scheme (\ref{2.16})  is unconditionally stable.
  \end{proof}

\subsection{The stability of the numerical methods in 3D}

\noindent{\bf Theorem 3.4.}
The D-AD scheme ({\ref{2.23})-({\ref{2.25}), defined by (\ref{2.22}), is unconditionally stable for   $\alpha,\beta,\gamma \in (1,2)$,
 if the matrices $A_{3D,{\Delta x}}^{k+1/2}$, $A_{3D,{\Delta y}}^{k+1/2}$ and $A_{3D,{\Delta z}}^{k+1/2}$ commute.

\begin{proof}

D-AD scheme (\ref{2.23})-(\ref{2.25}) can be written as
\begin{align}
&(I-A_{3D,{\Delta x}}^{k+1/2})U^{k,1}=(I+A_{3D,{\Delta x}}^{k+1/2}+2A_{3D,{\Delta y}}^{k+1/2}+2A_{3D,{\Delta z}}^{k+1/2})U^{k}+\Delta t F^{k+1/2}; \label{3.8}\\
&(I-A_{3D,{\Delta y}}^{k+1/2})U^{k,2}=U^{k,1}  - A_{3D,{\Delta y}}^{k+1/2}U^{k}; \label{3.9}\\
&(I-A_{3D,{\Delta z}}^{k+1/2})U^{k+1}=U^{k,2}  - A_{3D,{\Delta z}}^{k+1/2}U^{k}; \label{3.10}
\end{align}
according to (\ref{3.8})-(\ref{3.10}), we have the following equation
  \begin{equation}\label{3.11}
  \begin{split}
   &(I-A_{3D,{\Delta x}}^{k+1/2})(I-A_{3D,{\Delta y}}^{k+1/2})(I-A_{3D,{\Delta z}}^{k+1/2})U^{k+1} \\
  &\quad  =(I+A_{3D,{\Delta x}}^{k+1/2})(I+A_{3D,{\Delta y}}^{k+1/2})(I+A_{3D,{\Delta z}}^{k+1/2})U^{k}+\Delta tF^{k+1/2},
  \end{split}
  \end{equation}
  where the matrices $A_{3D,{\Delta x}}^{k+1/2}$ and $A_{3D,{\Delta y}}^{k+1/2}$ and $A_{3D,{\Delta z}}^{k+1/2}$ denote the
   operators $\frac{\Delta t}{2}\delta''_{\alpha,x}$ and $\frac{\Delta t}{2}\delta''_{\beta,y} $
    and $\frac{\Delta t}{2}\delta''_{\gamma,z} $, respectively, the vector $U^{k,1}$ and  $U^{k,2}$ denote the intermediate solution,
  the vector $F^{k+1/2}$ absorbs the source terms $f_{i,j,l}^{k+1/2}$ and the Dirichlet boundary conditions at time $t=t_{k+1}$ in the discretized equation. The matrices $A_{3D,{\Delta x}}^{k+1/2}$ and $A_{3D,{\Delta y}}^{k+1/2}$ and $A_{3D,{\Delta z}}^{k+1/2}$ are matrices of size $(N_x-1)(N_y-1)(N_z-1)\times(N_x-1)(N_y-1)(N_z-1)$.

 By the similar analysis, it can be proven that the spectral radius of the matrices $(I-A_{3D,{\Delta x}}^{k+1/2})^{-1}(I+A_{3D,{\Delta x}}^{k+1/2})$, $(I-A_{3D,{\Delta y}}^{k+1/2})^{-1}(I+A_{3D,{\Delta y}}^{k+1/2})$  and $(I-A_{3D,{\Delta z}}^{k+1/2})^{-1}(I+A_{3D,{\Delta z}}^{k+1/2})$ are less than 1,
therefore the difference  scheme (\ref{2.22})  is unconditionally stable.
\end{proof}

 \section{Multigrid method for the resulting matrix algebraic equations}
We use a V-cycle LOD-MGM to solve the resulting matrix algebraic equations of (\ref{1.1}).
Meanwhile, we show the  convergence  of the resulting system.
In order to develop a fast algorithm, i.e., realizing the computational count $\mathcal{O}(N \mbox{log} N)$ and the required storage $\mathcal{O}(N)$, as did in \cite{Pang:12, Wang:10, Wang:12}, we first introduce the Toeplitz matrix and the circulant matrix. The $n \times n$ Toeplitz  matrix $T_n(c)$ is defined by  \cite{Bottcher:05}
\begin{equation}\label{4.1}
T_n(c):=[c_{j-k}]_{j,k=1}^n=\left [ \begin{matrix}
                      c_0           &      c_{-1}             &      \cdots         &       c_{-(n-1)}       \\
                      c_{1}         &      c_{0}              &      \cdots         &       c_{-(n-2)}        \\
                     \vdots         &      \vdots             &      \ddots         &        \vdots            \\
                     c_{n-1}        &      c_{n-2}            &      \cdots         & c_{0}
 \end{matrix}
 \right ],
\end{equation}
and the circulant matrices are the ``periodic counsins" of Toeplitz matrices. We denote by $circ ~(c_0,c_1,\dots,c_{n-1})$ the circulant matrix
whose first column is $\widetilde{c}=(c_0,c_1,\dots,c_{n-1})^{\rm T}$,
\begin{equation}\label{4.2}
C_n:=\left [ \begin{matrix}
                      c_0           &      c_{n-1}      &      c_{n-2}       &      \cdots       &       c_2   &       c_1       \\
                      c_{1}         &      c_{0}        &      c_{n-1}       &      \cdots       &       c_3   &       c_{2}      \\
                      c_{2}         &      c_{1}        &      c_{0}         &      \ddots       &      \ddots &       c_{3}       \\
                      \vdots        &      \vdots       &      \ddots        &      \ddots       &      \ddots &       \vdots       \\
                      c_{n-2}       &      c_{n-3}      &      \ddots        &      \ddots       &       c_{0} &       c_{n-1}       \\
                      c_{n-1}       &      c_{n-2}      &      c_{n-3}       &      \cdots       &       c_{1} &       c_{0}
 \end{matrix}
 \right ].
\end{equation}
Moreover, we set $\omega_n=\mbox{exp}(2\pi i/n)$ and put
\begin{equation*}
F_n=\frac{1}{\sqrt{n}}\left [ \begin{array}{llllr}
                      1          &          1           &               1       &      \cdots         &               1             \\
                      1          &      \omega_n        &      \omega_n^{2}     &      \cdots         &       \omega_n^{n-1}         \\
                      1          &      \omega_n^{2}    &      \omega_n^{4}     &      \cdots         &       \omega_n^{2(n-1)}       \\
                      \vdots     &      \vdots          &      \vdots           &      \ddots         &        \vdots                  \\
                      1          &      \omega_n^{n-1}  &      \omega_n^{2(n-1)}&      \cdots         &        \omega_n^{(n-1)(n-1)}
\end{array}
 \right ],
\end{equation*}
with $i$ as the imaginary unit, and the matrix $F_n$ is called the Fourier matrix.
Therefore, a circulant matrix can be diagonalized by the  Fourier matrix $F_n$, i.e.,
\begin{equation}\label{4.3}
C_n:=F_n^{*}\mbox{diag}(F_n\widetilde{c} )F_n,
\end{equation}
where  $\mbox{diag}(F_n\widetilde{c} )$  is a diagonal matrix holding the eigenvalues of $C_n$.
From (\ref{4.3}), we can determine $\mbox{diag}(F_n\widetilde{c} )$
 in $\mathcal{O}(N \mbox{log} N)$ operations by the FFT of the first column $\widetilde{c} $ of $C_n$ \cite{Bottcher:05}.

 \subsection{ A $\mathcal{O}(N \mbox{log} N)$ V-cycle  MGM for 1D}
We employ the  V-cycle MGM  to solve
the one dimensional  system (\ref{2.11}) and illustrate the computational count of $\mathcal{O}(N \mbox{log} N)$ per iteration
and the required storage of $\mathcal{O}(N )$.

Suppose $A_h=I-A^{k+1/2}$, $u_h=U^{k+1}$ and $f_h=(I+A^{k+1/2})U^{k}+\Delta t F^{k+1/2}$,
then the  resulting system (\ref{2.11}) becomes the following  general linear system
\begin{equation}\label{4.4}
    A_hu_h=f_h,
  \end{equation}
and the  system (\ref{4.4})  can be carry out by the Algorithm \ref{V-cycle} \citep[p.\,443]{Saad:03} and \ref{MGM}.

In Algorithm \ref{V-cycle}, at the  highest (finest grid) level a mesh-size of $h$ is used to solve  the resulting  system (\ref{4.4}).
The finest grid operator  $A_h$ is with the finest grid size $h=\Delta x$;
the coarse grid operator   $A_H=I-A_{{2^lh}}$,
 $H=2^lh$, for $1\leq l \leq \mbox{log}_2N-1$;
   $h_0$ is the coarsest mesh-size;
   $I_H^h$, $I_h^H$  are respectively the prolongation operator
and the restriction operator.
For one dimensional system, the restriction operator$I_h^H$ is defined by \cite{Saad:03}
\begin{equation}\label{4.5}
I_h^H=\frac{1}{4}\left [ \begin{matrix}
                      1   &  2   &   1       &         &          &         &          &          & \\
                          &      &   1       &   2     &   1      &         &          &          &  \\
                          &      & \cdots    &  \cdots &  \cdots  &         &          &          &   \\
                          &      &           &         &          &     1   &    2     &     1    &    \\
 \end{matrix}
 \right ],
\end{equation}
and the prolongation operator $I_H^h=2(I_h^H)^{\rm T}$.  The smoothing operator $\mathtt{smooth}$ may be written as
\begin{equation}\label{4.6}
  \mathtt{smooth}(A_h,u_0,f_h)=S_hu_0+(I-S_h)A_h^{-1}f_h,
\end{equation}
where $S_h$ is the iteration matrix of the smoothing operator,
and we define the weighted (damped) Jacobi iteration matrix by \citep[p.\,9]{Briggs:00}
\begin{equation}\label{4.7}
  S_{h,\omega}=I-\omega D^{-1}A_h,
\end{equation}
with the weighting  factor $\omega \in \mathbb{R}$, and $D$ is the diagonal of $A_h$.
Thus, the (\ref{4.6}) becomes the following weighted  Jacobi iteration
\begin{equation}\label{4.8}
  \mathtt{smooth}(A_h,u_0,f_h)= S_{h,\omega}u_0+\omega D^{-1}f_h.
\end{equation}
In Algorithm \ref{V-cycle}, the factors $\nu_1$ and $\nu_2$ of $\mathtt{smooth}^{\nu_1}(A_h,u_0,f_h)$
 and $\mathtt{smooth}^{\nu_2}(A_h,u_h,f_h)$ denote the number of weighted  Jacobi iterations. In Algorithm \ref{MGM},  we give the stopping criterion of Algorithm \ref{V-cycle}.

Next, we illustrate the storage requirement of $\mathcal{O}(N)$ and the computational count of $\mathcal{O}(N \mbox{log} N)$ per iteration.

From (\ref{2.12}), we have $A^{k+1/2}={\rm diag}({\xi}^{k+1/2}) \widetilde{A}^{k+1/2}$, where

\begin{equation}\label{4.9}
\widetilde{A}^{k+1/2}=\left [ \begin{matrix}
2g_1^{\alpha}            &g_0^{\alpha}+g_2^{\alpha}&g_3^{\alpha}             &      \cdots              &  g_{N_x-2}^{\alpha}     &  g_{N_x-1}^{\alpha}  \\
g_0^{\alpha}+g_2^{\alpha}&        2g_1^{\alpha}    &g_0^{\alpha}+g_2^{\alpha}&       g_3^{\alpha}       &     \cdots              &  g_{N_x-2}^{\alpha} \\
g_3^{\alpha}             &g_0^{\alpha}+g_2^{\alpha}&2g_1^{\alpha}            & g_0^{\alpha}+g_2^{\alpha}&     \ddots              & \vdots  \\
\vdots                   &          \ddots         &       \ddots            &        \ddots            &      \ddots             &  g_3^{\alpha}  \\
g_{N_x-2}^{\alpha}       &          \ddots         &       \ddots            &        \ddots            &   2g_1^{\alpha}         & g_0^{\alpha}+g_2^{\alpha} \\
g_{N_x-1}^{\alpha}       &    g_{N_x-2}^{\alpha}   &    g_3^{\alpha}         &         \cdots           &g_0^{\alpha}+g_2^{\alpha}& 2g_1^{\alpha}
 \end{matrix}
 \right ]
\end{equation}
being a Toeplitz matrix, and $\xi^{k+1/2}=[\xi_1^{k+1/2},\xi_2^{k+1/2}\ldots,\xi_{N_x-1}^{k+1/2}]^{ T} $.

Then, we only need to store $\xi^{k+1/2}$
and $g^{\alpha}=[2g_1^{\alpha}, g_0^{\alpha}+g_2^{\alpha}, g_3^{\alpha} \ldots,g_{N_x-1}^{\alpha}]^{ T}   $
which have $2N-2$ parameters, instead of the full matrix $A^{k+1/2}$ which has $(N-1)^2$ parameters,
i.e., the required storage $\mathcal{O}(N)$.
Consider a one dimensional grid with $N$ points, the finest grid, $\Omega^h$, requires  $\mathcal{O}(N)$ storage locations;
$\Omega^{2h}$ requires $2^{-1}$ times as much storage as $\Omega^h$;  $\Omega^{4h}$ requires $4^{-1}$ times as much storage as $\Omega^h$;
in general, $\Omega^{ph}$ requires $p^{-1}$ times as much storage as $\Omega^h$.
Adding these terms we obtain \cite{Briggs:00}
\begin{equation*}
  \mbox{Storage}=\mathcal{O}(N) \cdot \left( 1+\frac{1}{2}+\frac{1}{2^2}+\ldots,+\frac{1}{2^{\log_2N-1}} \right)=\mathcal{O}(N).
\end{equation*}

Taking $v$ a given vector, for the Toeplitz matrix vector multiplication $\widetilde{A}^{k+1/2} v$,
 we first embed $\widetilde{A}^{k+1/2}$
into a $(2N_x-2) \times (2N_x-2) $ circulant matrix, i.e.,
 \begin{equation}\label{4.10}
\left [\begin{matrix}
                      \widetilde{A}^{k+1/2}          &     \ast             \\
                      \ast                              & \widetilde{A}^{k+1/2}
 \end{matrix} \right ]
 \left [\begin{matrix}
                     v                 \\
                     0
 \end{matrix}\right ]
 =\left [\begin{matrix}
                     \widetilde{A}^{k+1/2}v                \\
                     \dagger
 \end{matrix}\right ].
\end{equation}
Thus, using  (\ref{4.3}), the computational count of $ \widetilde{A}^{k+1/2}v$ remains as $\mathcal{O}(N \mbox{log} N)$.
Therefore, for a V-cycle MGM, each level is visited $\mathcal{O}(N \mbox{log} N)$ and grid $\Omega^{ph}$
requires $p^{-1}$ work units.  Similarly, adding these count we have
\begin{equation*}
\begin{split}
  &\mbox{V-cycle MGM computational count}\\
 &\quad  =\mathcal{O}(N \mbox{log} N) \cdot \left( 1+\frac{1}{2}+\frac{1}{2^2}+\ldots,+\frac{1}{2^{\log_2N-1}} \right)=\mathcal{O}(N \mbox{log} N).
  \end{split}
\end{equation*}

By similar analysis, we know that the required storage is still $\mathcal{O}(N)$ and the computational count $\mathcal{O}(N \mbox{log} N)$ for multidimensional case.

\subsection{ A $\mathcal{O}(N \mbox{log} N)$  V-cycle LOD-MGM for 2D}
 For D-AD scheme (\ref{3.3})-(\ref{3.4}), take
\begin{equation*}
\begin{split}
&A_{h_x}=I-A_{2D,{\Delta x}}^{k+1/2};  \quad u_{h_x}=U^{*}; \quad ~~~ f_{h_x}=(I+A_{2D,{\Delta x}}^{k+1/2}+2A_{2D,{\Delta y}}^{k+1/2})U^{k}+\Delta tF^{k+1/2}; \\
&A_{h_y}=I-A_{2D,{\Delta y}}^{k+1/2};  \quad  u_{h_y}=U^{k+1}; \quad  f_{h_y}=U^{*}-A_{2D,{\Delta y}}^{k+1/2}U^{k}.
\end{split}
\end{equation*}

Similarly, for PR-AD scheme (\ref{3.5})-(\ref{3.6}), denote
\begin{equation*}
\begin{split}
&A_{h_x}=I-A_{2D,{\Delta x}}^{k+1/2};  ~~ u_{h_x}=U^{*}; ~~~~~ f_{h_x}=(I+A_{2D,{\Delta y}}^{k+1/2})U^{k}+\frac{\Delta t}{2}F^{k+1/2}; \\
&A_{h_y}=I-A_{2D,{\Delta y}}^{k+1/2};  ~~ u_{h_y}=U^{k+1}; ~~ f_{h_y}=(I+A_{2D,{\Delta x}}^{k+1/2})U^{*}+\frac{\Delta t}{2}F^{k+1/2}.
\end{split}
\end{equation*}
Then, both the D-AD and PR-AD schemes,
defined by  (\ref{3.7}),  reduce to the following LOD form:
\begin{equation}\label{4.11}
    A_{h_x}u_{h_x}=f_{h_x},
  \end{equation}
  \begin{equation}\label{4.12}
    A_{h_y}u_{h_y}=f_{h_y}.
  \end{equation}

Therefore, we can solve the two dimensional system (\ref{3.7})
by V-cycle LOD-MGM (see Appendix Algorithm \ref{V-cycle}-\ref{LOD-MGM for 2D}).

The Algorithm \ref{LOD-MGM for 2D} starts with the initial time $t=0$ and  executes as follows:
  \begin{description}
 \item[(1)] First for every fixed  $y=y_j$ $ (j=1,\ldots,{N_y}-1)$, using Algorithm \ref{V-cycle} to solve a set of ${N_x}-1$ equations defined by (\ref{4.11}) at the mesh points $ x_i,i=1,\ldots,{N_x}-1$,
       to get $u_{h_x}$;
 \item[(2)] Next alternating the spatial direction, and for each fixed  $x=x_i$ $(i=1,\ldots,{N_x}-1)$ solving a set of ${N_y}-1$  equations defined by (\ref{4.12}) at the points $y_j,j=1,\ldots,{N_y}-1$,
  once  again we employ Algorithm \ref{V-cycle} to get $u_{h_y}$.
\end{description}


\subsection{ A $\mathcal{O}(N \mbox{log} N)$  V-cycle LOD MGM for 3D}
 For D-AD scheme (\ref{3.8})-(\ref{3.10}), similarly, we set
\begin{equation*}
\begin{split}
&A_{h_x}=I-A_{3D,{\Delta x}}^{k+1/2};  \quad u_{h_x}=U^{k,1}; \\
& f_{h_x}=(I+A_{3D,{\Delta x}}^{k+1/2}+2A_{3D,{\Delta y}}^{k+1/2}+2A_{3D,{\Delta z}}^{k+1/2})U^{k}+\Delta tF^{k+1/2}; \\
&A_{h_y}=I-A_{3D,{\Delta y}}^{k+1/2};  \quad  u_{h_y}=U^{k,2}; \quad  f_{h_y}=U^{k,1}-A_{3D,{\Delta y}}^{k+1/2}U^{k} ; \\
&A_{h_z}=I-A_{3D,{\Delta z}}^{k+1/2};  \quad  u_{h_z}=U^{k+1}; \quad  f_{h_z}=U^{k,2}-A_{3D,{\Delta z}}^{k+1/2}U^{k} ;
\end{split}
\end{equation*}

Then, the D-AD scheme
defined by (\ref{3.11}),  becomes the following form:
\begin{equation}\label{4.013}
    A_{h_x}u_{h_x}=f_{h_x},
  \end{equation}
  \begin{equation}\label{4.014}
    A_{h_y}u_{h_y}=f_{h_y},
  \end{equation}
  \begin{equation}\label{4.015}
    A_{h_z}u_{h_z}=f_{h_z}.
  \end{equation}
Therefore, we can   solve the three dimensional system (\ref{3.11})
by  Algorithm \ref{V-cycle}, \ref{MGM} and \ref{MGM for 3D}.

%

 \subsection{ Convergence analysis}
In this subsection, we discuss the convergence of  LOD-MGM.
For the convenience of convergence analysis, we assume the coefficient $c(x,y,z,t)$ is constant, and then all $\xi_i^{k+1/2}$ defined above are equal, and it can be denoted as $\xi$.
Let's first consider the one dimensional case.

From (\ref{4.9}), there exists
\begin{equation}\label{4.13}
A_h=I-A^{k+1/2}=I-\xi \widetilde{A}^{k+1/2} \equiv [c_{|j-k|}]_{(N_x-1)\times (N_x-1)},
\end{equation}
where
\begin{equation}\label{4.14}
  c_0=1- 2g_1^{\alpha} \xi; ~ c_1=-(g_0^{\alpha}+g_2^{\alpha}) \xi;~
 c_k=g_{k+1}^{\alpha} \xi,~~ k=2,\ldots,N_x-2,
\end{equation}
and $A_h$ is a symmetric Toeplitz matrix.
From the proof of Theorem 3.2, we know that the matrix $A_h$ is symmetric and strongly diagonally dominant
with positive diagonally elements, then $A_h$ is symmetric positive definite \citep[p.\,3]{Briggs:00}.

Since the matrix $A_h$ is symmetric positive define, we can define the following
three different  inner products \citep[p.\,78]{Ruge:87}
\begin{equation}\label{4.15}
  ( u,v  )_0=(Du,v), \quad (u,v)_1=(A_hu,v), \quad (u,v)_2=(D^{-1}A_hu,A_hv),
\end{equation}
where $D$ is the diagonal of $A_h$ and
along with their corresponding norms  $|| \cdot ||_i$ $(i=0,1,2)$.
If taking the coarsest grid size $H=h_0=2h$,
Algorithm  \ref{V-cycle} is called the two-grid method (TGM).
The TGM is rarely used in practice since the coarse grid operator may still be too
large to be solved exactly.  However, it is useful from a theoretical point view as the
first step to study the MGM convergence usually begins from the TGM \cite{Pang:12,Ruge:87,Saad:03, Wesseling:92}. Since the MGM convergence analysis is still a challenge topic in computational mathematics \cite{Briggs:00}.
In the following we only consider the convergence of the TGM.

\noindent{\bf Lemma 4.1.}
\emph{Let $A_h$ be a M-matrix and  the  weighted factor  $0 < \omega \leq 1 $, then weighted  Jacobi iteration (\ref{4.8}) converges.}
\begin{proof}
Taking $M=D/\omega$ and $N=D/\omega-A_h$,   since $0 < \omega \leq 1 $, then $M$ and $N$ be a regular splitting of a matrix $A_h$.
Note that $A_h$ is an M-matrix \citep[p.\,3]{Briggs:00},
thus we have \citep[p.\,119]{Saad:03}
\begin{equation*}
 \rho ( S_{h,\omega})=\rho (I-\omega D^{-1}A_h) < 1.
\end{equation*}
\end{proof}

\noindent{\bf Lemma 4.2} \citep[p.\,84]{Ruge:87}.
\emph{Let $A_h$ be a symmetric positive definite and  $\eta_0 \geq \rho(D^{-1}A_h)$.
If $\sigma \leq \omega(2-\omega \eta_0) $,
then the Jacobi relaxation with relaxation parameter $0<\omega <2/\eta_0$ satisfies
\begin{equation}\label{4.16}
 ||S_{h,\omega}e_h||_1^2 \leq  ||e_h||_1^2 - \sigma ||e_h||_2^2, \quad  \forall e_h \in \mathbb{R}^{N_x-1}.
\end{equation} }
The inequality (\ref{4.16}) is called the smoothing condition.  
For the TGM, the correction operator is given by \citep[p.\,85]{Briggs:00}
\begin{equation*}
  T_c=I-I_H^h(A_H)^{-1}I_h^HA_h,
\end{equation*}
therefore, the convergence factor of the TGM is $||(S_{h,\omega})^{\nu_2}\cdot T_{c}(S_{h,\omega})^{\nu_1}||_1$;
see \citep[p.\,89]{Ruge:87}. For convenience, we take $\nu_1=0$ and $\nu_2=1$. Therefore, the convergence factor
of the TGM is given by $||S_{h,\omega}\cdot T_c||_1$.

\noindent{\bf Lemma 4.3.} \citep[p.\,89]{Ruge:87}
\emph{Let $A_h$ be a symmetric positive definite matrix and  $S_{h,\omega}$ satisfy (\ref{4.16}).
Suppose that the interpolation $I_H^h$ has full rank and that, for each $e_h$,
  \begin{equation}\label{4.18}
   \min_{e^H \in \mathbb{R}^{N_x/2-1}}||e_h-I_H^he_H||_0^2\leq \kappa ||e_h||_1^2, \quad  \forall e_h \in \mathbb{R}^{N_x-1},
  \end{equation}
with  $\kappa>0$ independent of $e_h$. Then, $\kappa\geq \sigma$ and the convergence factor of the TGM convergence factor  satisfies
$||S_{h,\omega}\cdot T_c||_1 \leq \sqrt{1-\sigma/\kappa }$.}

Letting $L_{N_x-1}=\mbox{tridiag}(-1,2,-1)$ be the $(N_x-1)\times(N_x-1)$ one dimensional discrete laplacian, then $L_{N_x-1}$
is a symmetric positive definite matrix.
 We define $A_{rest}=A^{k+1/2}+c_1L_{N_x-1}$, where $c_1$ is defined by (\ref{4.14}) and it can also be shown  that $A_{rest}$
 is symmetric and diagonally dominant with positive  diagonally elements; then $A_{rest}$ is positive definite.
Hence, we have the following equation
\begin{equation}\label{4.19}
  (e_h,A_he_h)=(e_h,(I-c_1L_{N_x-1}+A_{rest})e_h) \geq  (e_h,(I-c_1L_{N_x-1})e_h),\quad  \forall e_h \in \mathbb{R}^{N_x-1}.
\end{equation}

\noindent{\bf Theorem 4.4.}
\emph{Since $A_h$, defined by (\ref{4.4}), is a symmetric positive definite matrix, if taking $\sigma=\omega(2-\omega\eta_0) $ with $\omega \in (0,1]$, then $S_{h,\omega}$ satisfies (\ref{4.16}) and the convergence factor of the TGM satisfies
\begin{equation}\label{4.20}
||S_{h,\omega}\cdot T_c||_1 < \sqrt{1-2\sigma/5 }<1,
\end{equation}}
where $\eta_0=\rho(  D^{-1}A_h)<2$.

\begin{proof}
From Lemma 4.1, we have
$
\rho(  D^{-1}A_h)   < 2.
$
Taking $\eta_0$ and $\sigma $ in Lemma 4.2 as $\rho(  D^{-1}A_h)$ and $\sigma= \omega(2-\omega\eta_0) $, respectively, then $S_{h,\omega}$ satisfies (\ref{4.16}).

Similar to the proof given in \cite{Pang:12}, we denote
 $e_h=(e_1,e_2,\ldots,e_{N_x-1})^{\rm T} \in \mathbb{R}^{N_x-1}$, $e_0=e_{N_x}=0$, and
$e_H=(e_2,e_4,\ldots,e_{N_x-2})^{\rm T} \in \mathbb{R}^{N_x/2-1}$;
 the norm $||\cdot||_0$ is defined by  (\ref{4.15}),
and $D=\mbox{diag}(A_h)=c_0I$,
where $c_0$ is defined by (\ref{4.14}).  There exists
 \begin{equation*}
||e_h-I_H^he_H||_0^2 = c_0\!\! \sum_{i=0}^{N_x/2-1}\!\!\!\left(e_{2i+1}-\frac{1}{2}e_{2i}-\frac{1}{2}e_{2i+2} \right)^2
                      \leq  c_0\!\! \sum_{i=1}^{N_x-1}\!\! \left(e_{i}^2 -e_{i}e_{i+1}\right).
\end{equation*}
From the above inequality, we obtain
\begin{equation*}
\begin{split}
\sum_{i=1}^{N_x-1}\!\! e_i^2  \geq \sum_{i=1}^{N_x-1}\!\! e_{i}e_{i+1}.\\
\end{split}
\end{equation*}
Similarly, we can check that
\begin{equation}\label{4.21}
\begin{split}
\sum_{i=1}^{N_x-1}\!\! e_i^2  \geq  -\sum_{i=1}^{N_x-1}\!\! e_{i}e_{i+1}.\\
\end{split}
\end{equation}
Combining (\ref{4.19}) with (\ref{4.21}), we have
\begin{equation*}
\begin{split}
||e_h||_1^2&=(e_h,A_he_h) \geq  (e_h,(I-c_1L_{N_x-1})e_h)=\sum_{i=1}^{N_x-1}\!\!\left( (1-2c_1)e_i^2+2c_1e_ie_{i+1}   \right)\\
&=\!\! \sum_{i=1}^{N_x-1}\!\!\left( \left(\frac{1}{2}-2c_1 \right) (e_i^2-e_ie_{i+1})+\frac{1}{2}(e_i^2 +e_ie_{i+1})   \right)
\geq \left(\frac{1}{2}-2c_1 \right)\!\! \sum_{i=1}^{N_x-1}\!\!\left(  (e_i^2-e_ie_{i+1})   \right).\\
\end{split}
\end{equation*}
Hence, Eq. (\ref{4.18}) holds,  i.e.,
\begin{equation*}
\begin{split}
||e_h-I_H^he_H||_0^2 \leq  \kappa ||e_h||_1^2,
\end{split}
\end{equation*}
\begin{equation*}
where
\begin{split}
\kappa=\frac{c_0}{\frac{1}{2}-2c_1}=\frac{2c_0}{1-4c_1} =\frac{2+4\xi(-g_1^{\alpha} )}{1+4\xi(g_0^{\alpha}+g_2^{\alpha})}
  \in \left(1,\frac{5}{2}\right),
\end{split}
\end{equation*}
since, according to Lemma 3.1, it is easy to check that $g_0^{\alpha}+g_2^{\alpha}< -g_1^{\alpha}<\frac{5}{2}(g_0^{\alpha}+g_2^{\alpha})$.
From Lemma 4.3,
we have  $||S_{h,\omega}\cdot T_c||_1 < \sqrt{1-2\sigma/5 }<1$.
\end{proof}
Theorem 4.4 shows that the V-cycle scheme has a convergence factor with bound,  that is independent of $\Delta x$.

For the multidimensional case, we can similarly define the weighted Jacobi iteration matrix
\begin{equation*}
  S_{h_x,\omega}=I-\omega D_{h_x}^{-1}A_{h_x}, \quad S_{h_y,\omega}=I-\omega D_{h_y}^{-1}A_{h_y}, \quad S_{h_z,\omega}=I-\omega D_{h_z}^{-1}A_{h_z}
\end{equation*}
and the LOD-TGM correction operators
\begin{equation*}
  T_{c_x}=I-I_H^h(A_H)^{-1}I_h^HA_{h_x}, \quad T_{c_y}=I-I_H^h(A_H)^{-1}I_h^HA_{h_y}, \quad T_{c_z}=I-I_H^h(A_H)^{-1}I_h^HA_{h_z},
\end{equation*}
then the convergent factors of the LOD-TGM satisfy
$ ||S_{h_x,\omega}\cdot T_{c_x}||_1 < \sqrt{1-2\sigma/5 }<1$, $||S_{h_y,\omega}\cdot T_{c_y}||_1 < \sqrt{1-2\sigma/5 }<1$ and $||S_{h_z,\omega}\cdot T_{c_z}||_1 < \sqrt{1-2\sigma/5 }<1$.

%
\noindent{\bf Remark 4.7 }  As mentioned in the Introduction section, nowadays there are two different second-order discretization schemes for fractional operators \cite{Sousa:12, Tian:12};  all the analysis given in this paper can be parallel extended to the case that the fractional operators are discretized by the scheme given in \cite{Tian:12}, in fact the scheme given in \cite{Tian:12} has more wide applications because of its good properties; and in Table \ref{tab:4}, we show the numerical results obtained by using the scheme of \cite{Tian:12} to discretize the fractional operators of (\ref{2.13}).


\section{Numerical results}
We employ the V-cycle MGM and V-cycle LOD-MGM described in Section 4 to solve the one dimensional case (\ref{2.8}) and multidimensional case (\ref{2.13},\ref{1.1}), respectively. The stopping criterion is taken as
\begin{equation*}
  \frac{||r^{(l)}||_2}{||r^{(0)}||_2}<10^{-7},
\end{equation*}
where $r^{(l)}$ is the residual vector after $l$ iterations.
In all tables, $N_t$ denotes the number of time steps; $N_x$, $N_y$ and  $N_z$, respectively, denotes the number of spatial grid points in $x$, $y$, and $z$ direction, and the numerical errors are measured by the $ l_{\infty}$ norm,  `Rate' denotes the convergent orders.
`CPU' denotes the total CPU time in seconds (s) or minutes (m) for solving the resulting discretized   systems,
 and `Iter'  denotes the average number of iterations required to solve a general linear system
$A_hu_h=f_h$  at each time level.

All numerical experiments are programmed in Python, and each computation was carried out  on a PC with the configuration: AMD Phenom (tm) II X4 830 CPU 2.79 GHZ and 3 GB RAM and a Linux operating system.  All the numerical results listed in the following tables are got by the V-cycle MGM or V-cycle LOD-MGM with the parameters:
the  number of  iterations $(\nu^1,\nu^2)=(1,1)$ and $(\omega_{pre},\omega_{post})=(1,1/2)$.

\noindent{\bf Remark 5.1.} From our numerical experiences, we find that:
 \begin{description}
 \item[(1)] With the increasing of the order of fractional derivative $\alpha$ from 1 to 2, the condition number of the matrix $\widetilde{A}^{k+1/2}$ becomes bigger and bigger; and when getting the same accuracy the cost for $\alpha=1.9$ almost double the cost for $\alpha=1.1$;
 \item[(2)] For making MGM more efficient, the parameters can be dynamically chosen as when $\alpha$ increases from 1 to 2, correspondingly $\omega_{post}$ decreases from 1 to 0.5 and fixing $\omega_{pre}=1$ or fixing $\omega_{post}=1$ and $\omega_{pre}$ decreases from 1 to 0.5;
 \item[(3)] MGM is still powerful for second-order schemes, when simulating (\ref{2.8}) with the parameters given in subsection 5.1, in Table \ref{tab:1} it is shown that the second-order scheme costs  $4.82\, s$  to obtain the accuracy with the maximum error $1.2407e-006$, and the first-order scheme used in \cite{Pang:12} costs $413.95\, s$ when getting the accuracy with the maximum error $2.0358e-006$.
 \end{description}


\subsection{Numerical results for 1D}
Let us consider  the one dimensional Riesz fractional  diffusion
equation (\ref{2.8}),
where $0< x < 1 $ and $0 < t \leq1$, with the variable coefficient $c(x,t)=x^{\alpha}t$, the forcing function
\begin{equation*}
\begin{split}
 f(x,t)&=-e^{-t}x^2(1-x)^2 \\
 &\quad +\frac{x^{\alpha}te^{-t}}{cos(\alpha \pi/2)}
   \left[ \frac{x^{2-\alpha}+(1-x)^{2-\alpha}}{\Gamma(3-\alpha)}   -6\frac{x^{3-\alpha}+(1-x)^{3-\alpha}}{\Gamma(4-\alpha)}
   +12\frac{x^{4-\alpha}+(1-x)^{4-\alpha}}{\Gamma(5-\alpha)} \right],
  \end{split}
\end{equation*}
and the initial condition $u(x,0)=x^2(1-x)^2 $ and the boundary
conditions $u(0,t)=u(1,t)=0$. This fractional PDE has the exact
value $u(x,t)=e^{-t}x^2(1-x)^2$, which may be confirmed by applying
the fractional differential equations
\begin{equation*}
\begin{split}
&_{x_L}D_{x}^{\nu}(x-x_L)^p=\frac{\Gamma(p+1)}{\Gamma(p+1-\nu)}(x-x_L)^{p-\nu},\\
&_xD_{x_R}^{\nu}(x_R-x)^p=\frac{\Gamma(p+1)}{\Gamma(p+1-\nu)}(x_R-x)^{p-\nu}.
\end{split}
\end{equation*}


\begin{table}[h]\fontsize{9.5pt}{12pt}\selectfont
  \begin{center}
  \caption{MGM to solve the resulting system (\ref{2.11}) of
  the 1D Riesz fractional convection diffusion equation (\ref{2.8}) at $t=1$ and $N_t=N_x$.}\vspace{5pt}

    \begin{tabular*}{\linewidth}{@{\extracolsep{\fill}}*{9}{c}}                                    \hline  
$N_t,\,N_x$ &  $\alpha=1.1$  & Rate           & Iter    & CPU   &  $\alpha=1.9$  &   Rate        &  Iter   &  CPU    \\\hline
    $2^5$   &   7.8755e-005  &     &  4.0     & 0.19 s             &  7.5578e-005   &          &  6.0    & 0.28 s\\\hline 
    $2^6$   &   2.1801e-005  & 1.8530    &  4.0     & 0.45 s       &  1.9255e-005   &  1.9727  &  6.0    & 0.74 s \\\hline 
    $2^7$   &   5.6999e-006  & 1.9354    &  4.0     & 1.29 s      &  4.8923e-006   &  1.9766  &  6.0    & 2.14 s \\\hline 
    $2^8$   &   1.4565e-006  & 1.9684    &  3.0     & 2.56 s       &  1.2407e-006   &  1.9794  &  6.0    & 4.82 s \\\hline 
    $2^9$   &   3.8540e-007  & 1.9181    &  3.0     & 7.45 s      &  3.1420e-007   &  1.9814  &  6.0    & 13.56 s \\\hline 
 $2^{10}$   &   9.7292e-008  & 1.9860    &  3.0     & 18.63 s     &  8.1028e-008   &  1.9552  &  6.0    & 34.99 s \\\hline 
    \end{tabular*}\label{tab:1}
  \end{center}
\end{table}


From Table \ref{tab:1}, we numerically confirm that the numerical scheme has second-order accuracy and the computational cost is of $\mathcal{O}(N \mbox{log} N)$ operations.


\subsection{Numerical results for 2D}
 Consider the two dimensional Riesz fractional convection diffusion equation (\ref{2.13}),
on a finite domain $ 0<  x< 1,\,0<  y< 1$, $0< t \leq 1$, and with the variable coefficients
 \begin{equation*}
\begin{split}
c(x,y,t)=x^{\alpha}y, \quad d(x,y,t)=xy^{\beta},
\end{split}
\end{equation*}
and the initial condition $u(x,y,0)=x^2(1-x)^2y^2(1-y)^2$ and the
Dirichlet boundary conditions on the rectangle in the form $
u(0,y,t)=u(x,0,t)=0$ and $ u(1,y,t)=u(x,1,t)=0$
 for all $ t\geq 0$. The exact solution to this two dimensional Riesz fractional convection diffusion equation is
 \begin{equation*}
\begin{split}
u(x,y,t)=e^{-t}x^2(1-x)^2y^2(1-y)^2.
  \end{split}
  \end{equation*}
From the above given quantities, it is easy to obtain the forcing function $f(x,y,t)$.

\begin{table}[h]\fontsize{9.5pt}{12pt}\selectfont
  \begin{center}
  \caption{LOD-MGM to solve the resulting system of
  the 2D Riesz fractional convection diffusion equation (\ref{2.13}) by the D-AD scheme  (\ref{2.17})-(\ref{2.18})  at $t=1$ and $N_t=N_x=N_y$.}\vspace{5pt}

    \begin{tabular*}{\linewidth}{@{\extracolsep{\fill}}*{9}{c}}                                    \hline  
$N_t,\,N_x,\,N_y$ &  $\alpha=1.1,\beta=1.1$  & Rate             & Iter    & CPU      &  $\alpha=1.8,\beta=1.9$  &   Rate        &  Iter   &  CPU    \\\hline
    $2^4$   &   2.4698e-005  &        &  4.5     & 2.09 s      &    2.5475e-005   &       &  7.0    & 3.51 s  \\\hline 
    $2^5$   &   6.1249e-006  & 2.0117   &  4.0     & 11.26  s     & 6.5211e-006   &  1.9659  &  6.0    & 18.03 s\\\hline 
    $2^6$   &   1.5212e-006  & 2.0095    &  4.0     & 55.85  s     & 1.6662e-006   &  1.9686  &  6.0    & 90.73 s\\\hline 
    $2^7$   &    3.7812e-007  & 2.0083    &  4.0     & 5 m 17  s    & 4.2362e-007    & 1.9757  & 6.0    & 8 m 45 s\\\hline 
    $2^8$   &    9.4076e-008  & 2.0070    &  3.0     & 21 m 36 s    &1.0744e-007   &  1.9792  &  6.0    & 39 m 1 s \\\hline 
    \end{tabular*}\label{tab:2}
  \end{center}
\end{table}


\begin{table}[h]\fontsize{9.5pt}{12pt}\selectfont
  \begin{center}
  \caption{LOD-MGM to solve the resulting system of
  the 2D Riesz fractional convection diffusion equation (\ref{2.13}) by the PR-AD scheme  (\ref{2.19})-(\ref{2.20}) at $t=1$ and $N_t=N_x=N_y$.}\vspace{5pt}

    \begin{tabular*}{\linewidth}{@{\extracolsep{\fill}}*{9}{c}}                                    \hline  
$N_t,\,N_x,\,N_y$ &  $\alpha=1.1,\beta=1.1$  & Rate             & Iter    & CPU       &  $\alpha=1.8,\beta=1.9$  &   Rate        &  Iter   &  CPU    \\\hline
    $2^4$   &   2.4698e-005   &          & 4.5     &   2.05 s          &  2.5475e-005   &            & 7.0      &  3.48 s      \\\hline 
    $2^5$   &    6.1249e-006  & 2.0117   & 4.0     &  11.17 s          &  6.5211e-006   &  1.9659    & 6.0      &  17.71 s     \\\hline 
    $2^6$   &   1.5212e-006   & 2.0095   & 4.0     &  55.14 s          &  1.6662e-006   &   1.9686   & 6.0      &  90.24 s     \\\hline 
    $2^7$   &   3.7812e-007   & 2.0083   & 4.0     &  5 m 18 s       &  4.2362e-007   &  1.9757    & 6.0      &  8 m 42 s    \\\hline 
    $2^8$   &   9.4075e-008   & 2.0070   & 4.0     &  22 m 25 s      &  1.0744e-007   &   1.9792   & 6.0      &  39 m 31 s   \\\hline 
    \end{tabular*}\label{tab:3}
  \end{center}
\end{table}



\begin{table}[h]\fontsize{9.5pt}{12pt}\selectfont
  \begin{center}
  \caption{LOD-MGM to solve the resulting system of the 2D Riesz fractional convection diffusion equation (\ref{2.13}) by the D-AD scheme  (\ref{2.17})-(\ref{2.18}) at $t=1$ and $N_t=N_x=N_y$; {\bf the Riesz fractional derivative (\ref{1.2}) is discretized by the scheme given in \cite{Tian:12}}.}\vspace{5pt}

    \begin{tabular*}{\linewidth}{@{\extracolsep{\fill}}*{9}{c}}                                    \hline  
$N_t,\,N_x,\,N_y$ &  $\alpha=1.1,\beta=1.1$  & Rate             & Iter    & CPU      &  $\alpha=1.8,\beta=1.9$  &   Rate        &  Iter   &  CPU    \\\hline
    $2^4$   &   2.4592e-005  &            &  5.0     & 2.25 s       &  2.4532e-005   &       &  7.0    & 3.5 s \\\hline 
    $2^5$   &   6.1745e-006  & 1.9938     &  4.0     & 11.42 s      &  6.0897e-006   &  2.0102  &  6.0    & 17.98 s \\\hline 
    $2^6$   &   1.5426e-006  & 2.0010     &  4.0     & 56.57 s      &  1.5102e-006   &  2.0116  &  6.0    & 90.74 s \\\hline 
    $2^7$   &   3.8444e-007  & 2.0045     &  4.0     & 5 m 28 s      &   3.7350e-007    &  2.0155  &  6.0    & 8 m 41 s \\\hline 
    $2^8$   &   9.5743e-008  & 2.0055     &  4.0     & 22 m 43 s    &  9.2374e-008   &  2.0155  &  6.0    & 38 m 22 s \\\hline 
    \end{tabular*}\label{tab:4}
  \end{center}
\end{table}

 From Table \ref{tab:2} and Table \ref{tab:3}, numerically it can also be noticed that the D-AD and PR-AD are equivalent in two dimensional problems.  We employ the LOD-MGM to solve  two dimensional Riesz fractional diffusion equation, numerical results further display the computational cost is of $\mathcal{O}(N \mbox{log} N)$ operations and the numerical scheme is second-order convergent.

In particular, we further numerically confirm that this paper is still valid if the Riesz fractional derivative (\ref{1.2}) is discretized by another existing second-order discretization scheme given in \cite{Tian:12}, see Table \ref{tab:4}. In fact, theoretically we can also easily draw the same conclusion.


\subsection{Numerical results for 3D}
Consider the three dimensional Riesz fractional convection diffusion equation (\ref{1.1}),
on a finite domain $ 0<  x< 1$, $0<  y< 1$, $0<  z< 1$, and $0< t \leq 1$, and with the variable coefficients
 \begin{equation*}
\begin{split}
c(x,y,z,t)=x^{\alpha}yz, \quad d(x,y,z,t)=xy^{\beta}z, \quad e(x,y,z,t)=xyz^{\gamma},
\end{split}
\end{equation*}
and the initial condition $u(x,y,z,0)=x^2(1-x)^2y^2(1-y)^2z^2(1-z)^2$ and the
zero Dirichlet boundary conditions on the cube. The exact solution to this three dimensional Riesz fractional convection diffusion equation is
 \begin{equation*}
\begin{split}
u(x,y,z,t)=e^{-t}x^2(1-x)^2y^2(1-y)^2z^2(1-z)^2.
  \end{split}
  \end{equation*}
According to the  above conditions, it is easy to obtain the forcing function $f(x,y,z,t)$.


\begin{table}[h]\fontsize{9.5pt}{12pt}\selectfont
  \begin{center}
  \caption{LOD-MGM to solve the scheme (\ref{2.23})-(\ref{2.25}) of
  the 3D Riesz fractional convection diffusion equation (\ref{1.1}) at $t=1$ and $N_t=N_x=N_y=N_z$.}\vspace{5pt}

    \begin{tabular*}{\linewidth}{@{\extracolsep{\fill}}*{9}{c}}                                    \hline  
$N_t$ &  $\alpha=\beta=\gamma=1.1$  & Rate             & Iter    & CPU      &  $\alpha=1.8,\beta=1.9,\gamma=1.8$  &   Rate        &  Iter   &  CPU \\\hline
    $2^3$   &   5.9349e-006   &           &  4.75     &3.31 s      & 5.8311e-006   &        &  7.0    &5.28 s \\\hline 
    $2^4$   &  1.4792e-006  & 2.0044    &  4.0     & 42.19   s    &  1.4867e-006   &  1.9717  & 7.0    & 70.1 s\\\hline 
    $2^5$   &   3.7377e-007  & 1.9846    &  4.5     & 7 m 40 s       &  3.8428e-007   &  1.9519  &  6.0    &13 m 1 s \\\hline 
    $2^6$   &   9.3376e-008  & 2.0010    &  4.37     & 76 m 13 s       &  9.8179e-008   &  1.9687 &  6.0    &136 m 1 s \\\hline 
    \end{tabular*}\label{tab:5}
  \end{center}
\end{table}


The numerical results in Table \ref{tab:5} are obtained by employing the LOD-MGM to solve the three dimensional Riesz fractional diffusion equation, they again display the  computational count of $\mathcal{O}(N \mbox{log} N)$ operations
and the scheme is second-order convergent.

%

\section{Conclusions}

With the appearing of the two effective second-order discretization schemes \cite{Sousa:12, Tian:12} and MGM being successfully employ to solve the resulting system of the one dimensional fractional diffusion equation discretized by first-order scheme \cite{Pang:12}, our attentions turn to the possibility of efficiently solving the multidimensional fractional Riesz diffusion equation by second-order scheme and MGM. This paper shows that when solving (\ref{1.1}) the second-order schemes are unconditionally stable, and the structure of the resulting matrix algebraic equations is almost the same as the one by the first-order scheme and then none computational costs increase but accuracy is greatly improved if using the second scheme instead of the first-order one. And LOD-MGM still preserves its powerfulness of  $\mathcal{O}(N \mbox{log} N)$ computational counts and of $\mathcal{O}(N)$ storage when solving the resulting matrix system of the multidimensional fractional Riesz diffusion equation discretized by the second-order scheme. For computing the completely same equation, from Table \ref{tab:1}, it can be noticed that the second-order scheme costs  $4.82\, s$  to obtain the accuracy with the maximum error $1.2407e-006$, and the first-order scheme used in \cite{Pang:12} costs $413.95\, s$ when getting the accuracy with the maximum error $2.0358e-006$.

Last but not least, we want to refer to that although this paper focus on using the second-order discretization given in \cite{Sousa:12}, all the analysis of this paper is still valid if applying the second-order discretization in \cite{Tian:12}; in fact the validness is already verified numerically in Table \ref{tab:4}.

\section*{Acknowledgments}This work was supported by the Program for
New Century Excellent Talents in University under Grant No.
NCET-09-0438, the National Natural Science Foundation of China under
Grant No. 10801067 and No. 11271173, and the Fundamental Research Funds for the
Central Universities under Grant No. lzujbky-2010-63 and No. lzujbky-2012-k26.

\section*{Appendix}
\label{sec:V-cycle}
For a general linear system
\begin{equation*}
A_hu_h=f_h,
\end{equation*}
we  employ the following
 V-cycle MGM   (Algorithm \ref{V-cycle}-\ref{MGM}) to solve one dimensional (\ref{2.8})
and V-cycle LOD-MGM (Algorithm \ref{V-cycle}-\ref{LOD-MGM for 2D}) solve two dimensional (\ref{2.13}).
Solve the three-dimensional  system (\ref{1.1})
by  Algorithm \ref{V-cycle},\ref{MGM} and \ref{MGM for 3D}.
\begin{algorithm}
\caption{MGM for 1D ~~~~$u_h$=V-cycle$(A_h,u_0,f_h)$ }
\label{V-cycle}
\begin{algorithmic}[1]
\STATE Pre-smooth:\quad $u_h:=\mathtt{smooth}^{\nu_1}(A_h,u_0,f_h)$
\STATE Get residual:\quad $r_h=f_h-A_hu_h$
\STATE Coarsen:\quad $r_H=I_h^Hr_h$
\IF{$H==h_0$}
\STATE Solve:\quad $A_H\xi_H=r_H$
\ELSE
\STATE Recursion:\quad $\xi_H=\mbox{V-cycle}(A_H,0,r_H)$
\ENDIF
\STATE Correct:\quad $u_h:=u_h+I_H^h\xi_H$
\STATE Post-smooth:\quad $u_h:=\mathtt{smooth}^{\nu_2}(A_h,u_h,f_h)$
\STATE Return $u_h$
\end{algorithmic}
\end{algorithm}

\begin{algorithm}
\caption{Stopping criterion for MGM}
\label{MGM}
\begin{algorithmic}[1]
\STATE $u_h:=u_0$
\STATE $r0:=||A_hu_0-f_h||_2$
\STATE $r_h:=r_0$
\WHILE{$\frac{r_h}{r_0}>\epsilon$}
\STATE $u_h:=\mbox{V-cycle}(A_h,u_h,f_H)$
\STATE $r_h:=||A_hu_h-f_h||_2$
\ENDWHILE
\STATE Return $u_h$
\end{algorithmic}
\end{algorithm}

\begin{algorithm}[t]
\caption {LOD-MGM for 2D}
\label {LOD-MGM for 2D}
\begin{algorithmic}[1]
\STATE $t:=0$
\WHILE{$t<T$}
\STATE $t:=t+\Delta t$
\FOR{$\mbox{every fixed}\  y_j,(j=1:N_y-1)$}
\STATE solve system (\ref{4.11}) by Algorithm \ref{MGM}
\ENDFOR
\FOR{each fixed $x_i,(i=1:N_x-1)$}
\STATE solve system (\ref{4.12}) by Algorithm \ref{MGM}
\ENDFOR
\ENDWHILE
\end{algorithmic}
\end{algorithm}

\begin{algorithm}[t]
\caption {LOD-MGM for 3D}
\label {MGM for 3D}
\begin{algorithmic}[1]
\STATE $t:=0$
\WHILE{$t<T$}
\STATE $t:=t+\Delta t$
\FOR{$\mbox{every fixed}\  z_l,(l=1:N_z-1)$ and  $y_j,(j=1:N_y-1)$}
\STATE solve system (\ref{4.013}) by Algorithm \ref{MGM}
\ENDFOR
\FOR{each fixed $x_i,(i=1:N_x-1)$ and $z_l,(l=1:N_z-1)$}
\STATE solve system (\ref{4.014}) by Algorithm \ref{MGM}
\ENDFOR
\FOR{each fixed $y_j,(j=1:N_y-1)$ and $x_i,(i=1:N_x-1)$}
\STATE solve system (\ref{4.015}) by Algorithm \ref{MGM}
\ENDFOR
\ENDWHILE
\end{algorithmic}
\end{algorithm}

\section*{}

\end{document}